\documentstyle[11pt]{article}
\newcommand\nc{\newcommand} \nc\rnc{\renewcommand} 

\nc\nn{\newenvironment} \nc\nt{\newtheorem}
\nt{Claim}{Claim}[section]
\nt{claim}{Claim}[subsection]
\nt{Corollary}{Corollary}[section]
\nt{corollary}{Corollary}[subsection]
\nt{DefinitionAux}{Definition}[section]
\nn{Definition}{\begin{DefinitionAux}\em}{\ $\Box$\end{DefinitionAux}}
\nt{definitionaux}{Definition}[subsection]
\nn{definition}{\begin{definitionaux}\em}{\ $\Box$\end{definitionaux}}
\nt{Lemma}{Lemma}[section]
\nt{lemma}{Lemma}[subsection]
\nt{ProofAux}{Proof} 
\nn{Proof}{\begin{ProofAux}\em}{\ $\Box$\end{ProofAux}}
\nt{Theorem}{Theorem}[section]
\nt{theorem}{Theorem}[subsection]

\nc\w{\wedge} \nc\st{\,:\:}

\def\E {\mathop{\mbox{\lower.3ex\hbox{\large E}}}}	
\def\EE{\mathop{\mbox{\lower.7ex\hbox{\LARGE  E}}}}	
\rnc\P{\mbox{\bf P}}			

\rnc\a{\alpha} \rnc\b{\beta} \nc\g{\gamma} \nc\e{\varepsilon} 
\nc\X{\bar{X}} \nc\Y{\bar{Y}}
\nc\LL{L^\omega_{\infty,\omega}} \nc\Lk{L^k_{\infty,\omega}}

\begin{document}

\title{On Finite Rigid Structures}

\author{Yuri Gurevich\thanks{Partially supported by BSF, NSF and ONR.
Electrical Engineering and Computer Science, University of Michigan, Ann
Arbor, MI 48109-2122, USA} \and Saharon Shelah\thanks{Partially supported
by BSF and NSF.  Mathematics, Hebrew University, Jerusalem 91904, Israel,
and Mathematics, Rutgers University, New Brunswick, NJ 08903, USA}}
\date{}\maketitle

\begin{abstract}
The main result of this paper is a probabilistic construction of finite
rigid structures.  It yields a finitely axiomatizable class of finite rigid
structures where no $\LL$ formula with counting quantifiers defines a
linear order.  
\end{abstract}

\section{Introduction}

In this paper, structures are finite and of course vocabularies are finite
as well.  A class is always a collection of structures of the same
vocabulary which is closed under isomorphisms.  

An $r$-ary {\em global relation\/} on a class $K$ is a function $\rho$ that
associates an $r$-ary relation $\rho_A$ with each structure $A\in K$ in
such a way that every isomorphism from $A$ to a structure $B$ extends to an
isomorphism from the structure $(A,\rho_A)$ to the structure $(B,\rho_B)$
[G].

Recall that a structure is {\em rigid\/} if it has no nontrivial
automorphisms.  If a binary global relation $<$ defines a linear order in a
class $K$ (that is, on each structure in $K$) then every structure in $K$
is rigid.  Indeed, suppose that $\theta$ is an automorphism of a structure
$A\in K$ and let $a$ be an arbitrary element of $A$.  Since
\begin{eqnarray*}
A \models \theta(x)<\theta(a) &\iff& A \models x<a, \\
A \models \theta(x)>\theta(a) &\iff& A \models x>a, 
\end{eqnarray*} 
the number of elements preceding $\theta(a)$ in the linear order $<_A$
equals the number of elements preceding $a$.  Hence $\theta(a)=a$.

Conversely, if every structure in a class $K$ is rigid then some global
relation $\rho$ defines a linear order on each structure in $K$.  Alex
Stolboushkin constructed a finitely axiomatizable class of rigid structures
such that no first-order formula defines a linear order in $K$ [S].  Anuj
Dawar conjectured that, for every finitely axiomatizable class $K$ of rigid
structures, some formula in the fixed-point extension of first-logic
defines a linear order in $K$ [D].  Using the probabilistic method, we
refute the conjecture and construct a finitely axiomatizable class of
structures where no $\LL$ formula with counting quantifiers defines a
linear order (Theorem~4.1).  At the end of Section~4, we answer a question
of Scott Weinstein [W] related to rigid structure.

To make this paper self-contained, we provide a reminder in the rest of
this section.  As in a popular version of first-order logic, $\LL$ formulas
are built from atomic formulas by means of negations, conjunctions,
disjunctions, the existential quantifier and the universal quantifier.  The
only difference is that, in $\LL$, one is allowed to form the conjunction
and the disjunction of an arbitrary set $S$ of formulas provided that the
total number of variables in all $S$-formulas is finite.  $\LL(C)$ is the
extension of $\LL$ by means of {\em counting quantifiers\/} $(\exists2x)$,
$(\exists3x)$, etc.  The semantics is obvious.  $\Lk$ (resp.\ $\Lk(C)$) is
the fragment of $\LL$ (resp.\ $\LL(C)$) where formulas use at most $k$
variables.

There is a pebble game $G^k(A,B)$ appropriate to $\Lk(C)$ [IL].  Here $A$
and $B$ are structures of the same purely relational vocabulary.  The game
is played by Spoiler and Duplicator on a board comprised by $A$ and $B$.
For each $i=1,\ldots,k$, there are two identical pebbles marked by $i$.
Initially there are no pebbles on the board.  After every round, either
both $i$-pebbles are off the board or else one of them covers an element of
$A$ and the other covers an element of $b$; furthermore the pebbles on the
board define a partial isomorphism from $A$ to $B$.  (This means that
(i)~an $i$-pebble and a $j$-pebble cover different elements of $A$ if and
only if their twins cover different elements of $B$, and (ii)~the map that
takes a pebble-covered element of $A$ to the element of $B$ covered by the
pebble of the same number is a partial isomorphism.)

A round of $G^k(A,B)$ is played as follows.  
\begin{enumerate}
\item
Spoiler chooses a number $i$; if the $i$-pebbles are on the board, they are
taken off the board.  Then Spoiler chooses a structure $M\in\{A,B\}$ and a
nonempty subset $X$ of $M$.
\item
Duplicator chooses a subset $Y$ of the remaining structure $N$ such that
$\|Y\|=\|X\|$.  If $N$ has no subsets of cardinality $\|X\|$, the game is
over; Spoiler has won and Duplicator has lost.
\item
Spoiler puts an $i$-pebble on an element $y\in Y$.
\item
Duplicator puts the other $i$-pebble on an element $x\in X$ in such a way
that the pebbles define a partial isomorphism.  If $X$ has no appropriate
element $x$, the game is over; Spoiler has won and Duplicator has lost.
Otherwise Duplicator wins the round
\end{enumerate}

Spoiler wins a play of the game if the number of rounds in the play is
infinite.  

\begin{Theorem}[{[IL]}]
If Duplicator has a winning strategy in $G^k(A,B)$ then no $\Lk(C)$
sentence $\phi$ distinguishes between $A$ and $B$.
\end{Theorem}

It is not hard to prove the theorem by induction on $\phi$.  The converse
implication is true too [IL] but we will not use it.

\paragraph{Acknowledgment}
This investigation has been provoked by a stimulating conversation that one
of us had with Steven Lindell and Scott Weinstein at the beginning of
October~1993.

\section{Hypergraphs}

\subsection{Preliminaries}

In this paper, a {\em hypergraph\/} is a pair $H=(U,T)$ where $U=|H|$ is a
nonempty set and $T$ is a collection of 3-element subsets of $U$; elements
of $U$ are {\em vertices\/} of $H$, and elements of $T$ are {\em
hyperedges\/} of $H$.  It can be seen as a structure with universe $U$ and
irreflexive symmetric ternary relation $\{(x,y,z) \st \{x,y,z\}\in T \}$.

Every nonempty subset $X$ of $U$ gives a {\em sub-hypergraph\/} 
\[ H|X = ( X, \{ h\st h\in T \w h\subseteq X \} \]
of $H$.  The number of hyperedges in $H|X$ will be called the {\em
weight\/} of $X$ and denoted $[X]$.  As usual, the number of vertices of
$X$ is called the cardinality of $X$ and denoted $\|X\|$.

Vertices $x,y$ of a hypergraph $H$ are {\em adjacent\/} if there is a
hyperedge $\{x,y,z\}$; the vertex $z$ {\em witnesses\/} that $x$ and $y$
are adjacent.

\begin{definition}
A vertex set $X$ is {\em dense\/} if $\|X\|\leq2[X]$.  A hypergraph is {\em
$l$-meager\/} if it has no dense vertex sets of cardinality $\leq2l$.
\end{definition}

\begin{lemma}
In a 2-meager hypergraph, the intersection of any two distinct hyperedges
contains at most one vertex.
\end{lemma}

\begin{Proof}
If $\|h_1\cap h_2\|=2$ then $h_1\cup h_2$ is $2$-dense. 
\end{Proof}

\begin{definition}
A vertex set $X$ is {\em super-dense\/} or {\em immodest\/} if
$\|X\|<2[X]$.  A hypergraph is {\em $l$-modest\/} if it has no super-dense
sets of cardinality $\leq2l$.
\end{definition}

It follows that if $X$ is a dense vertex set of cardinality $\leq2l$ in an
$l$-modest hypergraph then $\|X\|=2[X]$ and in particular $\|X\|$ is even.

\subsection{Cycles}

\begin{definition}
A sequence $x_1,\ldots,x_k$ of $k\geq3$ distinct vertices is a {\em weak
cycle\/} of length $k$ if it satisfies the following two conditions where
the subscripts are viewed as numbers modulo $k$:
\begin{enumerate}\item
Each $x_i$ is adjacent to $x_{i+1}$.
\item
Either $k>3$ or else $k=3$ but $\{x_1,x_2,x_3\}$ is not a hyperedge.
\end{enumerate}
\end{definition}

We will index elements of a weak cycle of length $k$ with numbers modulo
$k$.

\begin{definition}
A weak cycle $x_1,\ldots,x_k$ is a {\em cycle\/} of length $k\geq3$ if no
triple $x_i,x_{i+1},x_{i+2}$ forms a hyperedge.  A corresponding {\em
witnessed cycle\/} of length $k$ is a vertex sequence
$x_1,\ldots,x_k,y_1,\ldots,y_k$ where each $y_i$ witnesses that $x_i$ is
adjacent to $x_{i+1}$.
\end{definition}

\begin{definition}
A vertex sequence $x_1,x_2$ is a {\em cycle\/} of length $2$ if there are
distinct vertices $y_1,y_2$ different from $x_1,x_2$ such that
$\{x_1,x_2,y_1\}$ and $\{x_2,x_1,y_2\}$ are hyperedges; the sequence
$x_1,x_2,y_1,y_2$ is a corresponding {\em witnessed cycle\/} of length $2$.
\end{definition}

\begin{lemma}
Every weak cycle includes a cycle.  More exactly, some (not necessarily
contiguous) subsequence of a weak cycle is a cycle.  Thus, an acyclic
hypergraph (that is, a hypergraph without any cycles) has no weak cycles.
\end{lemma}

\begin{Proof}
We prove the lemma by induction on the length.  Let $x_1,\ldots,x_k$ be a
weak cycle that is not a cycle, so that some $x_i,x_{i+1},x_{i+2}$ is a
hyperedge; without loss of generality, $i=1$.  Then the sequence
$x_1,x_3,\ldots,x_k$ of length $k-1$ is a weak cycle or a hyperedge.  In
the first case, use the induction hypothesis.  In the second, $k=4$ and
$x_1,x_3$ form a cycle witnessed by $x_2$ and $x_4$.
\end{Proof}

\begin{theorem}
In any $l$-modest graph, 
\begin{itemize}\item 
every minimal dense set of cardinality $2k\leq2l$ is a witnessed cycle of
length $k$, and
\item
every witnessed cycle of length $k\leq l$ is a minimal dense set of
cardinality $2k$.
\end{itemize}
\end{theorem}

The theorem clarifies the structure of minimal dense sets of cardinality
$\leq2l$ which play an important role in our probabilistic construction.
However the theorem itself will not be used and can be skipped.  The rest
of this subsection is devoted to proving the theorem.

\begin{Proof}
Fix some number $l\geq2$ and restrict attention to $l$-modest hypergraphs.

\begin{lemma}
For every vertex set $X$, the following statements are equivalent:
\begin{enumerate}
\item $X$ is a dense set of cardinality $4$.
\item $X$ is a minimal dense set of cardinality $4$
\item Vertices of $X$ form a witnessed cycle of length $2$. 
\end{enumerate}
\end{lemma}

\begin{Proof}
It is easy to see that (1) is equivalent to (2) and that (3) implies (1).
It remains to check that (1) implies (3).  Suppose (1).  By $l$-modesty
$[X]=2$.  Thus, $X$ includes two hyperedges $h_1$ and $h_2$.  Clearly,
$h_1\cup h_2=X$ and $\|h_1\cap h_2\|=2$.  It is easy to see that the
vertices of $h_1\cap h_2$ form a cycle and the vertices of $X$ form a
corresponding witnessed cycle.
\end{Proof}

In the rest of this subsection, $3\leq k\leq l$.

\begin{lemma}
Every witnessed cycle $x_1,\ldots,x_k,y_1,\ldots,y_k$ forms a dense set of
cardinality $2k$.
\end{lemma}

\begin{Proof}
Let $W=\{x_1,\ldots,x_k,y_1,\ldots,y_k\}$.  It suffices to check that the
$k$ hyperedges $\{x_i,x_{i+1},y_i\}$ are all distinct.  For then, using
$l$-modesty, we have
\[ 2k \leq 2[W] \leq \|W\| \leq 2k. \]

If $i\neq j$ but $\{x_i,x_{i+1},y_i\} =\{x_j,x_{j+1},y_j\}$ then either
$x_j=x_{i+1}$ or else $x_j=y_i$ in which case $x_{j+1}=x_i$.  Without loss
of generality, $x_j=x_{i+1}$ and therefore $j=i+1$ modulo $k$.  If also
$x_{j+1}=x_i$ then $i=j+1=i+2$ modulo $k$ which contradicts the fact that
$k>2$.  Thus $x_{j+1}=y_i$, so that $y_i=x_{i+2}$ and therefore
$\{x_i,x_{i+1},x_{i+2}\}$ is a hyperedge which contradicts the definition
of cycles.
\end{Proof}

\begin{lemma}
Every minimal dense vertex set of cardinality $2k$ forms a witnessed cycle of
length $k$.
\end{lemma}

\begin{Proof}
Without loss of generality, the given minimal vertex set contains all
vertices of the given hypergraph $H$; if not, restrict attention to the
corresponding sub-hypergraph of $H$.

It suffices to prove that $H$ includes a weak cycle of length $\leq k$.
For then, by Lemma~2.2.1, $H$ includes a cycle of length $\leq k$.  If a
witnessed version of the cycle contains less than $2k$ vertices then, by
the previous lemma, $H$ contains a proper dense subset.  

By contradiction suppose that $H$ does not include a weak cycle of length
$k$.

\begin{claim}
A hypergraph of cardinality $2k$ is acyclic if no proper vertex set is
dense and there is no weak cycles of length $\leq k$.
\end{claim}

\begin{Proof}
By contradiction suppose that there is a cycle of length $m>k$ and choose
the minimal possible $m$.  Consider a witnessed cycle $x_1,\ldots,x_m,
y_1,\ldots,y_m$.

Since the hypergraph has $<2m$ vertices, some $y_i$ occurs in
$x_1,\ldots,x_m$.  Without loss of generality, $y_1=x_j$ for some $j$, so
that $\{x_1,x_2,x_j\}$ is a hyperedge and therefore $j$ differs from $1,2$
and $3$.  But then the sequence $x_2,\ldots,x_j$ is a weak cycle and thus
includes a cycle of length $<m$.  This contradicts the choice of $m$.
\end{Proof}

\begin{claim}
Any acyclic hypergraph of positive weight contains a hyperedge $Y$ such
that at most one vertex of $Y$ belongs to any other hyperedge.
\end{claim}

\begin{Proof}
Let $s=(x_1,\ldots,x_k)$ be a longest vertex sequence such that (i)~for
every $i<k$, $x_i$ is adjacent to $x_{i+1}$, and (ii)~for no $i<k-1$, the
triple $x_i,x_{i+1},x_{i+2}$ forms a hyperedge.  Since the hypergraph has
hyperedges, $k\geq2$.  If $k=2$ then all hyperedges are disjoint and the
claim is obvious.  Suppose that $k\geq3$.

Pick a vertex $y$ such that $Y=\{x_{k-1},x_k,y\}$ is a hyperedge.  We prove
that neither $x_k$ nor $y$ belongs to any other hypergraph.  Since there
are no cycles of length $2$, $y$ is uniquely defined.  We prove that
neither $x_k$ nor $y$ belongs to any other hypergraph.  Vertex $y$ does not
occur in $x_1,\ldots,x_k$; otherwise $x_i,\ldots,x_{k_1}$ is a weak cycle.
Notice that $y$ can replace $x_k$ in $s$.  Thus it suffices to prove that
$x_k$ does not belong to any other hyperedge.

By contradiction, suppose that a hyperedge $Z\neq Y$ contains $x_k$ and let
$z\in Z-Y$.  By the maximality of $s$, it contains $z$; otherwise $s$ can
be extended by $z$.  But then the final segment $S=[z,x_k]$ of $s$ forms a
weak cycle.
\end{Proof}

\begin{claim}
No acyclic hypergraph is dense.
\end{claim}

\begin{Proof}
Induction on the cardinality of the given hypergraph $I$.  The claim is
trivial if $[I]=0$.  Suppose that $[I]>0$.  By the previous claim, $I$ has
a hyperedge $X=\{x,y,z\}$ such that neither $y$ nor $z$ belongs to any
other hyperedge.  Let $J$ be the sub-hypergraph of $I$ obtained by removing
vertices $y$ and $z$.  Using the induction hypothesis, we have
\[ \|I\| = \|J\| + 2 > 2[J] + 2 = 2 ([I] + 1) = 2 [I]. \]
\end{Proof}

Now we are ready to prove the lemma.  By Claim~2.2.1, $H$ is acyclic.  By
Claim~2.2.3, $H$ is not dense which gives the desired contradiction.
\end{Proof}

\begin{lemma}
Every witnessed cycle of length $k$ forms a minimal dense set.
\end{lemma}

\begin{Proof}
Let $W$ be the set of the vertices of the given witnessed cycle of length
$k$.  By Lemma~2.2.3, $W$ is a dense set of cardinality $2k$.  By the
$l$-modesty of the hypergraph, $W$ contains precisely $k$ hyperedges.  It
is easy to see now that every proper subset $X$ of $W$ is acyclic; by
Claim~2.2.3, $X$ is not dense.
\end{Proof}

Lemmas~2.2.2--2.2.5 imply the theorem.
\end{Proof}

\subsection{Green and Red Vertices}

Fix $l\geq2$ and consider a sufficiently modest hypergraph.  More
precisely, we require that the hypergraph is $(2l+2)$-modest.  It follows
that, for every dense set $V$ of cardinality $\leq4l+4$, $\|V\|=2[V]$.

For brevity, we use the following terminology.  A minimal dense vertex set
of cardinality $\leq2l$ is a {\em red block\/}.  A vertex is {\em red\/} if
it belongs to a red block; otherwise it is {\em green\/}.  A hyperedge is
{\em green\/} if it consists of green vertices.  The {\em green
sub-hypergraph\/} is the sub-hypergraph of green vertices.

\begin{lemma}
Distinct red blocks are disjoint.
\end{lemma}

\begin{Proof}
We suppose that distinct red blocks $X$ and $Y$ have a nonempty
intersection $Z$ and prove that the union $V = X\cup Y$ is immodest.
Indeed, $Z$ is a proper subset of $X$; otherwise $Y$ is not a minimal dense
set.  Therefore $Z$ is not dense and
\[ \|V\| = \|X\| + \|Y\| - \|Z\| = 2[X] + 2[Y] - \|Z\|
< 2[X] + 2[Y] - 2[Z] = 2( [X] + [Y] - [Z]) \leq 2[V]. \]
\end{Proof}

\begin{lemma}
Adjacent red vertices belong to the same red block. 
\end{lemma}

\begin{Proof}
Suppose that adjacent red vertices $x$ and $y$ belong to different red
blocks $X$ and $Y$ respectively, and let $h$ be a hyperedge containing $x$
and $y$.  We show that the set $V=X\cup Y\cup h$ is immodest.  Indeed,
\[ \|V\| \leq \|X\| + \|Y\| + 1 = 2[X] + 2[Y] + 1
< 2([X] + [Y] + 1) \leq 2[V]. \]
\end{Proof}

\begin{lemma}
No green vertex is adjacent to two different red vertices.
\end{lemma}

\begin{Proof}
By contradiction suppose that a green vertex $b$ is adjacent to distinct
red vertices $x$ and $x'$.  Let $X,X'$ be the red blocks of $x,x'$
respectively, $h$ be a hyperedge containing $b$ and $x$, and $h'$ be a
hyperedge containing $b$ and $x'$.  We show that the set $V = X\cup X'\cup
h\cup h'$ is immodest.  By the previous lemma, $h=h'$ implies $X=X'$.

If $h=h'$ then
\[ \|V\| = \|X\|+1 = 2[X]+1 < 2([X]+1) \leq[V]. \]

If $h\neq h'$ but $X=X'$ then
\[ \|V\| \leq \|X\|+3 = 2[X]+3 < 2 ([X] + 2) \leq 2 [V]. \]

If $X\neq X'$ then
\[ \|V\| \leq \|X\| + \|X'\| + 3 = 2[X] + 2[X'] + 3
< 2 ([X] + [X'] + 2) \leq [V]. \]
\end{Proof}

\begin{definition}
A hypergraph is {\em odd\/} if, for every nonempty vertex set $X$, there is
a hyperedge $h$ such that $\|h\cap X\|$ is odd.  
\end{definition}

For future reference, some assumptions are made explicit in the following
theorem. 

\begin{theorem}
Suppose that a hypergraph $H$ of cardinality $n$ satisfies the following
conditions where $n'<n$.
\begin{itemize}
\item $H$ is $(2l+2)$-modest.
\item The number of red vertices is $<n'$.
\item Every vertex set of cardinality $\geq n'$ includes a hyperedge.
\item 
For every nonempty vertex set $X$ of cardinality $<n'$, there exist a
vertex $x\in X$ and distinct hyperedges $h_1, h_2$ such that $h_1\cap X =
h_2\cap X = \{x\}$.
\end{itemize}
Then the green sub-hypergraph of $H$ is an odd, $l$-meager hypergraph of
cardinality $>n-n'$.
\end{theorem}

\begin{Proof}
Since the green sub-hypergraph $G$ is obtained from $H$ by removing all
dense vertex sets of cardinality $\leq2l$, $G$ is $l$-meager.  By the
second condition, $\|G\|>n-n'$.  To check that $G$ is odd, let $X$ be a
nonempty set of green vertices.  If $\|X\|\geq n'$, use the third
condition.  Suppose that $\|X\|<n'$ and let $x,h_1,h_2$ be as in the fourth
condition; both $\|h_1\cap X\|$ and $\|h_2\cap X\|$ are odd.  By
Lemma~2.3.3, at least one of the two hyperedges is green.
\end{Proof}

\subsection{Attraction}

\begin{definition}
In an arbitrary hypergraph, a vertex set $X$ {\em attracts\/} a vertex $y$
if there are vertices $x_1, x_2$ in $X$ such that $\{x_1,x_2,y\}$ is a
hyperedge.  $X$ is {\em closed\/} if it contains all elements attracted by
$X$.  As usual, the {\em closure\/} $\X$ of $X$ is the least closed set
containing $X$.
\end{definition}

\begin{lemma}
In an $l$-meager hypergraph, if $X$ is a vertex set of cardinality $k\leq
l$ then $\|\X\|<2k$.
\end{lemma}

\begin{Proof}
Construct sets $X_0,\ldots,X_m$ as follows.  Set $X_0=X$.  Suppose that
sets $X_0,\ldots,X_i$ have been constructed.  If $X_i$ is closed, set $m=i$
and terminate the construction process.  Otherwise pick a hyperedge $h$
such that $\|h\cap X_i\|=2$ and let $X_{i+1}=h\cup X_i$.  We show that
$m<k$.  

By contradiction suppose that $m\geq k$.  Check by induction on $i$ that
$\|X_i\|=k+i$ and $[X_i]\geq i$.  Since the hypergraph is $l$-meager, we
have: $2[X_k] <\|X_k\| = 2k\leq 2[X_k]$.  This gives the desired
contradiction.
\end{Proof}

\begin{lemma}
Suppose that $Y$ is a vertex set of cardinality $\leq k$ in a $2k$-meager
hypergraph and $p=\|\Y-Y\|$.  Then $p<n$ and there is an ordering
$z_1,\ldots,z_p$ of $\Y-Y$ such that each $z_j$ is attracted by
$Y\cup\{z_i\st i<j\}$.
\end{lemma}

\begin{Proof}
By the previous lemma, $\|\Y\|<2\|Y\|$.  Hence $p = \|\Y-Y\| < \|Y\| \leq
n$.  Choose elements $z_j$ by induction on $j$.  Suppose that $1\leq j\leq
p$ and all elements $z_i$ with $i<j$ have been chosen.  Since $\|\Y\|
=\|\|Y\|\|+p$ vertices, the set $Z_{j-1} = Y\cup\{z_i\st i<j\}$ is not
closed.  Let $z_j$ be any element in $\Y-Y$ attracted by $Z_{j-1}$.
\end{Proof}

\begin{theorem}
Suppose that $X$ is a vertex set of cardinality $<k$ in a $2k$-meager
hypergraph, $z_0\notin\X$, $Y=\X\cup\{z_0\}$, $Z=\Y$ and $p=\|Z-Y\|$.  Then
$p<k$ and there is an ordering $z_1,\ldots,z_p$ of $Z-Y$ such that, for
every $j>0$, $z_j$ is attracted by $Y\cup\{z_i\st 1\leq i<j\}$ and there is
a unique hyperedge $h_j$ witnessing the attraction.
\end{theorem}

\begin{Proof}
By the previous lemma, $p<k$.  Construct sequence $z_1,\ldots,z_p$ as in
the proof of the previous lemma.  For any $j>0$, let $h_j$ be a hyperedge
witnessing that $Z_{j-1}=Y\cup\{z_i\st 1\leq i<j\}$ attracts $y_j$.

By contradiction suppose that, for some positive $j\leq p$, some hyperedge
$h_j'\neq h_j$ witnesses that $z_j$ is attracted by $Z_{j-1}$.  Let
$S=\{h_1,\ldots,h_j,h_j'\}$.  We show that $V=\bigcup S$ is a dense set of
cardinality $\leq2k$ which contradicts the $2k$-meagerness of the
hypergraph.  

Since $V$ contains all hyperedges in $S$, $[V]\geq j+1$.  Since none of the
vertices $z_1,\ldots,z_j$ is attracted by $\X$, $\|h\cap \X\|\leq1$ for all
$h\in S$ and thus $\|V\cap\X\|\leq j+1$.  We have
\[ \|V\| = \|(V\cap\X)\cup\{z_0,\ldots,z_j\}\| 
\leq (j+1)+(j+1) \leq 2\cdot[V]. \]
Thus $V$ is a dense set of cardinality $\|V\|\leq2(j+1)\leq2(p+1)\leq2k$.
\end{Proof}

\section{Existence}

\begin{Theorem}
For any integers $l\geq2$ and $N>0$, there exists an odd $l$-meager
hypergraph of cardinality $>N$.
\end{Theorem}

In fact, there exists an odd $l$-meager hypergraph of cardinality precisely
$N$ but we do not need the stronger result here. 

\begin{Proof}
Now fix $l\geq2$ and $N>0$ and choose a positive real $\e<1/(2l+3)$.  Let
$n$ range over integers $\geq2N$ divisible by $4$ and $U$ be the set of
positive integers $\leq n$.  For each 3-element subset $a$ of $U$, flip a
coin with probability $p=n^{-2+\e}$ of heads, and let $T$ is the collection
of triples $a$ such that the coin comes up heads.  This gives a random
graph $H=(U,T)$.

We will need the following simple inequality.  In this section,
$\exp\a=e^\a$ and $\log\a=\log_e\a$.

\begin{Claim}
For all positive reals $q,r,s$ such that $p^r<1/2$, 
\begin{eqnarray}
&& \exp(-2qn^{s-2r+r\e}) < (1-p^r)^{qn^s} < \exp(-qn^{s-2r+r\e})
\end{eqnarray}
\end{Claim}

\begin{Proof}
Suppose that $0<\a<1/2$.  By Mean Value Theorem applied to function
$f(t)=-\log(1-t)$ on the interval $[0,\a]$, there is a point $t\in(0,\a)$
such
\[ f(\a)-f(0) = -\log(1-\a) = (\a-0) f'(t) = \a/(1-t). \]
Since $\a < \a/(1-t) < \a/(1-\a) < a/(1-1/2) = 2\a$, we have $\a <
-\log(1-\a) < 2\a$ and therefore $e^{-2\a}<1-\a<e^{-\a}$.  Now let $\a=p^r$
and raise the terms to power $qn^s$.
\end{Proof}

Call an event $E=E(n)$ {\em almost sure\/} if the probability $\P[E]$ tends
to 1 as $n$ grows to infinity.  We prove that, almost surely, $H$ satisfies
the conditions of Theorem~2.3.1 with $n'=n/4$ and therefore the green
subgraph of $H$ is an odd $l$-meager graph of cardinality $>N$.

\begin{Lemma}
Almost surely, $H$ is $(2l+2)$-modest.
\end{Lemma}

\begin{Proof}
It suffices to prove that, for each particular $m\leq4l+4$, the probability
$q_m$ that there is a super-dense vertex sets of cardinality $m$ is $o(1)$.
A vertex set $X$ of cardinality $m$ is super-dense if $m<2[X]$, that is, if
$X$ includes more than $m/2$ hyperedges.  Let $k$ be the least integer that
exceeds $m/2$.  Then $m\leq2k-1$ and therefore $n^{m-2k}\leq n^{-1}$.  Also
$2k-2\leq m\leq4l+4$, so that $k\leq2l+3$ and $k\e<1$.  Let $M={m\choose3}$
and $c={M\choose k}$.  We have
\[ q_m  < {n\choose m} \cdot c\cdot p^k <
c\cdot n^m\cdot n^{(-2+\e)k} = c\cdot n^{m-2k+k\e} \leq
c\cdot n^{-1+k\e} = o(1). \]
\end{Proof}

\begin{Lemma}
Almost surely, the number of red vertices is $<n/4$.
\end{Lemma}

\begin{Proof}
It suffices to prove that the expected number of red vertices is $o(n)$.
Indeed, let $r$ be the number of red vertices and $s$ ranges over the
integer interval $[n/4,n]$.  Then
\[ \EE[r] \geq \sum s\cdot \P[r=s] 
\geq \frac{n}{4} \sum \P[r=s] = \frac{n}{4} \P[r\geq\frac{n}{4}] \]
and thus $\P[r\geq\frac{n}{4}]$ tends to $0$ if $\E[r]=o(n)$.

Furthermore, it suffices to show that, for each particular $m\leq2l$, the
expected number $f(m)$ of vertices $v$ such that $v$ belongs to a dense set
$X$ of cardinality $m$ is $o(n)$.  Let $k=\lceil m/2\rceil$.  Then
$m\leq2k$ and therefore $n^{m-2k}\leq1$.  Also, $2k\leq m-1<2l$ and
therefore $k<l$ and $k\e<1$.  Let $M={m\choose3}$ and $c={M\choose k}$.  We
have
\[ f(m) \leq n \cdot {n-1\choose m-1} c p^k < n\cdot n^{m-1} c p^k =
c \cdot n^m p^k = c\cdot n^{m-2k+k\e} \leq c\cdot n^{k\e} = o(n). \]
\end{Proof}

\begin{Lemma}
Almost surely, every vertex set of cardinality $\geq n/4$ includes a
hyperedge.
\end{Lemma}

\begin{Proof}
Chose a real $c>0$ so small that $cn^3\leq{n/4\choose3}$ and let $q$ be the
probability that there exists a vertex set of cardinality $\geq n/4$ which
does not include any hyperedges.  Using inequality (1), we have
\[ q < 2^n \cdot (1-p)^{n/4\choose3} < e^n \cdot (1-p)^{cn^3}
< e^n \cdot \exp(-cn^{1+\e}) = o(1). \]
\end{Proof}

\begin{Lemma}
For every nonempty vertex set $X$ of cardinality $<n/4$, there exist a
vertex $x\in X$ and hyperedges $h_1, h_2$ such that 
\[ h_1\cap X = h_2\cap X = h_1\cap h_2 = \{x\}. \]
\end{Lemma}

\begin{Proof}
Let $X$ range over nonempty vertex sets of cardinality $<n/4$, $Y$ be the
collection of even numbers $y\in U-X$, and $Z$ be the collection of odd
numbers $z\in U-X$.  Clearly, $\|Y\|\geq n/4$ and $\|Z\|\geq n/4$.

Let $x$ range over $X$, $\sigma(x,X)$ mean that there exist vertices
$y_1,y_2\in Y$ such that $\{x,y_1,y_2\}$ is a hyperedge, and $\tau(x,X)$
mean that there exist vertices $z_1,z_2\in Z$ such that $\{x,z_1,z_2\}$ is
a hyperedge.  Call $X$ {\em bad\/} if and $\sigma(x,X)\w\tau(x,X)$ fails
for all $x$.  We prove that, almost surely, there are no bad vertex sets.

Choose a real $c>0$ so small that $cn^2<{n/4\choose2}$.  For given $X$
and $x$, 
\[ \P[\neg\sigma(x,X)] = (1-p)^{\|Y\|\choose2} 
\leq (1-p)^{n/4\choose2} < (1-p)^{cn^2} < \exp[-cn^\e]. \] 
The last inequality follows from inequality (1).  Similarly,
$\P[\neg\tau(x,X)] <\exp[-cn^\e]$.  Hence
\[ \P[\neg\sigma(x,X)\vee\neg\tau(x,Y)] 
\leq \P[\neg\sigma(x,X)] + \P[\neg\tau(x,X)] < 2\exp[-cn^\e] =
\exp[\log2-cn^\e]. \]

If $\|X\|=m$ then
\[\P[X \mbox{ is bad }] < 
\big(\exp[\log2-cn^\e]\big)^m = \exp[m(\log2-cn^\e)]. \]

For each $m<n/4$, let $q_m$ be the probability that there is a bad vertex
set of cardinality $m$.  For sufficiently large $n$, $\log2n-cn^\e<0$ and
therefore $\exp(\log2n-cn^\e)<1$.  Thus
\[ q_m \leq n^m \cdot \exp[m(\log2-cn^\e)] =
\exp[ m (\log2n - cn^\e) ] \leq \exp[ \log2n - cn^\e ]. \]

Finally, let $q$ be the probability of the existence of a bad set.  We have
\[ q < \frac{n}{4} \exp[\log2n - cn^\e] = o(1). \]
\end{Proof}

Theorem~3.1 is proved.
\end{Proof}

\section{Multipedes}

The domain $\{x\st\exists y (xEy)\}$ and the range $\{y\st\exists x
(xEy)\}$ of a binary relation $E$ will be denoted $D(E)$ and $R(E)$
respectively.

\begin{Definition}
A {\em $1$-multipede\/} is a directed graph $(U,E)$ such that $D(E)\cap
R(E) =\emptyset$, $D(E)\cup R(E) = U$, every element in $D(E)$ has exactly
one outgoing edge and every element in $R(E)$ has exactly two incoming
edges.
\end{Definition}

If $xEy$ holds then $x$ is a {\em foot\/} of $y$ and $y$ is the {\em
segment\/} $S(x)$ of $x$.  We extend function $S$ as follows.  If $x$ is a
segment then $S(x)=x$.  If $X$ is a set of segments and feet then $S(X) =\{
S(x)\st x\in X\}$.

\begin{Definition}
A {\em $2^-$-multipede\/} is a structure $(U,E,T)$ such that $(U,E)$ is a
$1$-multipede and $(U,T)$ is a hypergraph where each hyperedge $h$
satisfies the following conditions:
\begin{itemize}\item
Either all elements of $h$ are segments or else all elements of $h$ are
feet.
\item If $h$ is a foot hyperedge then $S(h)$ is a hyperedge as well.
\end{itemize}
\end{Definition}

If $X=\{x,y,z\}$ is a segment hyperedge then every 3-element foot set $A$
with $S(A)=X$ is a {\em slave\/} of $X$.  A slave $A$ of $X$ is {\em
positive\/} if $A$ is a hyperedge; otherwise it is {\em negative\/}.  Two
slaves of $X$ are {\em equivalent\/} if they are identical or one can be
obtained from the other by permuting the feet of two segments.  In other
words, if $a,a'$ are different feet of $x$ and $b,b'$ are different feet of
$y$ and $c,c'$ are different feet of $z$ then the eight slaves of $X$ split
into the following two equivalence classes
\[ \{a,b,c\}, \{a,b',c'\}, \{a',b,c'\}, \{a',b',c\} \]
and
\[ \{a',b,c\}, \{a,b',c\}, \{a,b,c'\}, \{a',b',c'\} \]

\begin{Definition}
A {\em $2$-multipede\/} is a $2^-$-multipede where, for each segment
hyperedge $X$, exactly four slaves of $X$ are positive and all four
positive slaves are equivalent.
\end{Definition}

A $2$-multipede $(U,E,T)$ is {\em odd\/} if the segment hypergraph
$(R(E),T)$ is so.

\begin{Lemma}
If an automorphism $\theta$ of an odd $2$-multipede does not move any
segment then it does not move any foot either.
\end{Lemma}

\begin{Proof}
By contradiction suppose that $\theta$ moves a foot $a$ of a segment $x$.
Clearly, $\theta(a)$ is the other foot of $x$.  Let $X$ be the collection
of segments $x$ such that $\theta$ permutes the feet of $x$.  Since the
multipede is odd, there exists a segment hyperedge $h$ such that $\|h\cap
X\|$ is odd.  It is easy to see that $\theta$ takes positive slaves of $X$
to negative ones and thus is not an automorphism.
\end{Proof}

\begin{Lemma}
Let $M$ is a $2k$-meager $2$-multipede and $\Upsilon$ be the extension of
the vocabulary of $M$ by means of individual constants for every segment of
$M$.  No $\Lk(C)$ sentence in the vocabulary $\Upsilon$ distinguishes
between $M$ and the $2$-multipede $N$ obtained from $M$ by permuting the
feet of one segment.
\end{Lemma}

To be on the safe side, let us explain what it means that $N$ is obtained
from $M$ by permuting the feet of one segment.  To obtain $N$, choose a
segment $x$ and perform the following transformation for every segment
hyperedge $h$ that contains $x$: Make all positive slaves of $h$ negative
and the other way round.

\begin{Proof}
Call a collection  $X$ of segments and feet {\em closed\/} if it satisfies
the following conditions:
\begin{itemize}
\item
The segments of $X$ form a closed set in the sense of Definition~2.4.1.
\item
If $a$ is foot of $x$ then $a\in X\leftrightarrow x\in X$.
\end{itemize}

Call a partial isomorphism $\a$ from $M$ to $N$ {\em regular\/} if $\a$
leaves segments intact and takes any foot to a foot of the same segment.
The domain of a partial isomorphism $\a$ will be denoted $D(\a)$.  A
regular partial isomorphism $\a$ is {\em safe\/} if there is a regular
extension of $\a$ to the closure $\overline{D(\a)}$.

\begin{Claim}
Each safe partial isomorphism $\a$ from $M$ to $N$ has a unique regular
extension to $\overline{D(\a)}$.
\end{Claim}

\begin{Proof}
Let $X=D(\a)$ and suppose that $\b$ and $\g$ are regular extension of $\a$
to $\X$.  Let $Y=S(X)$ and $Z=S(\Y)$.  By Lemma~2.4.2, there exists a
linear order $z_1,\ldots,z_p$ of the elements of $Z-Y$ such that each $z_j$
is attracted by the set $Z_{i-1}=Y\cap\{y_i\st i<j\}$.  We need to prove
that, for every $j$, either both $\b$ and $\g$ leave the feet of $z_j$
intact or else both of them permute the feet.  We proceed by induction on
$j$.  Suppose that $\b$ and $\g$ coincide on the feet of every $y_i$ with
$i<j$ and let $h$ witness that $Z_{j-1}$ attracts $z_j$.  Let $\{a,b,c\}$
be any positive slave of $h$ where $c$ is a feet of $z_j$.  By the
induction hypothesis, $\b(a)=\g(a)$ and $\b(b)=\g(b)$; let $a'=\b(a)$ and
$b'=\b(b)$.  Since $\b$ and $\g$ are partial isomorphisms, both
$\{a',b',\b(c)\}$ and $\{a',b',\g(c)\}$ are hyperedges in $N$.  Since $N$
is a $2$-multipede, $\b(c)=\g(c)$.
\end{Proof}

The unique regular extension of $\a$ will be denoted $\bar{\a}$.

\begin{Claim}
Suppose that $\a$ is a safe partial isomorphism from $M$ to $N$ with domain
$X$ of cardinality $<n$.  For every element $a\in|M|-\X$, there is a safe
extension of $\a$ to $X\cup\{a\}$ which leaves $a$ intact.
\end{Claim}

\begin{Proof}
We construct a regular extension $\b$ of $\bar{\a}$ to
$\overline{X\cup\{a\}}$.  Let $z_0$ be the segment of $a$,
$Y=S(\X)\cup{z_0}$, $Z=S(\Y)$ and $p=\|Z-Y\|$.  By Theorem~2.4.1, there is
a linear ordering $z_1,\ldots,z_p$ on the vertices of $Z-Y$ such that, for
every $j>0$, $z_j$ is attracted by $Y\cup\{z_i\st 1\leq i<j\}$ and there is
a unique hyperedge $h_j$ witnessing the attraction.

The desired $\b$ leaves intact all segments in $Z$ and the feet of $z_0$.
It remains to define $\b$ on the feet of segments $z_j$, $1\leq j\leq k$.
We do that by induction on $j$.  Suppose that $\b$ is defined on the feet
of all $z_i$ with $i<j$ and let $h_j$ be as above.  Let $d$ be a foot of
$y_j$ and pick a positive slave $\{b,c,d\}$ of $h_j$ in $M$; $\b$ is
already defined at $b$ and $c$.  The slave $\{\b(b),\b(c),\b(d)$ of $h_j$
should be positive in $N$.  This defines uniquely whether $\b(d)$ equals
$d$ or the other foot of $y_j$.

We need to check that $\b$ is a partial isomorphism from $M$ to $N$.  The
only nontrivial part is to check that if $A$ is a slave of a segment
hyperedge $h$ then $A$ is positive in $M$ if and only if $\b(A)$ is
positive in $N$.  Without loss of generality, $A\not\subseteq\X$.  Let $j$
be the least number such that $S(\X)\cup\{z_0,\ldots,z_j\}$ includes $h$.
Since $\X$ does not attract $z_0$, $\X$ includes all hyperedges in
$S(\X)\cup\{z_0\}$; thus $j>0$.  By the uniqueness property of $h_j$,
$h=h_j$.  By the construction of $\b$, $A$ is positive in $M$ if and only
if $\b(A)$ is positive in $N$.
\end{Proof}

The desired winning strategy of Duplicator is to ensure that, after each
round, pebbles define a safe partial isomorphism.  Suppose that pebbles
define a safe partial isomorphism $\a$ and Spoiler starts a new round.  By
the symmetry between $M$ and $N$, we may suppose that Spoiler chooses $M$
and a subset $X$ of elements of $M$.  Duplicator chooses $N$ and a subset
$\{f(x)\st x\in X\}$ where $f$ is as follows.  If $x\in\overline{D(\a)}$
then $f(x)=\bar{\a}(x)$; otherwise $f(x)=x$.  Now use the previous Lemma.
\end{Proof}

\begin{Definition}
A {\em $3$-multipede\/} is a structure $(M,<)$ where $M$ is a $2$-multipede
and $<$ is a linear order on the set of segments of $M$.
\end{Definition}

\begin{Definition}
A {\em $4$-multipede\/} is a $3$-multipede together with (i)~additional
elements representing uniquely all sets of segments and (ii)~the
corresponding containment relation $\e$.  
\end{Definition}

We skip the details of the definition of $4$-multipedes.  The additional
elements are called {\em super-segments\/}.  

A $4$-multipede is {\em odd\/} if the hypergraph of segments is so.

\begin{Lemma}
The collection of odd $4$-multipedes is finitely axiomatizable.
\end{Lemma}

\begin{Proof}
We give only three axioms which express that every set of segments is
represented by a unique super-segment:
\begin{itemize}\item
There is a super-segment $Y$ such that there is no $x$ with $x\e Y$.
\item
For every super-segment $Y$ and every segment $x$, there exists a
super-segment $Y'$ such that, for every $y$, $y\e Y' \leftrightarrow (y\e Y
\vee y=x)$.
\item
Super-segments $Y$ and $Y'$ are equal if $x\e Y\leftrightarrow x\e Y'$ for
all $x$.
\end{itemize}
\end{Proof}

\begin{Lemma}
Every odd $4$-multipede is rigid.
\end{Lemma}

\begin{Proof}
Let $\theta$ is an automorphism of a $4$-multipede $M$.  Because of the
linear order on segments, $\theta$ leaves intact all segments.  Therefore
it leaves intact all super-segments.  By Lemma~4.1, it leaves intact all
feet as well.
\end{Proof}

A $4$-multipede is {$l$-meager\/} if the hypergraph of segments is so.

\begin{Lemma}
Let $M$ is a $2k$-meager $4$-multipede and $\Upsilon$ be the extension of
the vocabulary of $M$ by means of individual constants for every segment of
$M$.  No $\Lk(C)$ sentence in the vocabulary $\Upsilon$ distinguishes
between $M$ and the $4$-multipede $N$ obtained from $M$ by permuting the
feet of a segment.
\end{Lemma}

\begin{Proof}
The proof is similar to that of Theorem~4.1.  We use the terminology and
notation of the proof of Theorem~4.1.  Call a collection of segments, feet
and super-segments {\em closed\/} if the subcollection of segments and feet
is so.  Lemma~4.2 remains true.  Lemma~4.3 remains true as well; if $a$ is
a super-segment, then $\X\cup\{a\}$ is closed and the desired $\b$ is the
extension of $\bar{\a}$ by means of $\g(a)=a$.  The remainder of the proof
is as above.
\end{Proof}

\begin{Lemma}
There exists $j$ such that no $\Lk(C)$ formula defines a linear order in
any $2(j+k)$-meager $4$-multipede.
\end{Lemma}

\begin{Proof}
Let $M$ be any structure in the vocabulary of $4$-multipedes, $M'$ be an
extension of $M$ with individual constants for all elements of $M$, and
$N=M''$ be an extension of $M'$ with a linear order $<$.  There exists an
$\LL$ sentence $\psi_N$ which describes $N$ up to isomorphism: For each
basic relation $R$ of $N$ and each tuple $\bar{x}$ of elements of $M$ of
appropriate length, $\psi_N$ says whether $\bar{x}$ belongs to $R$ or not.
{\em Cf.\/} [HKL].  The number $j$ of variables in $\psi$ does not depend
on $M$.

By contradiction suppose that an $\Lk(C)$ formula $\phi$ defines a linear
order in an $2(j+k)$-meager $4$-multipede $M$.  Define $M'$ as above and
let $M$ be the extension of $M'$ by means of the linear order $<$ defined
by $\phi$.  Replace each atomic formula $t_1<t_2$ in $\psi_N$ with
$\phi(t_1,t_2)$; here each $t_i$ is a variable or an individual constant.
The resulting $L^{j+k}_{\infty,\omega}(C)$ formula describes $M'$ up to
isomorphism.  This contradicts the preceding lemma.
\end{Proof}

\begin{Theorem}
There exists a finitely axiomatizable class of rigid structures such that
no $\LL(C)$ sentence that defines a linear order in every structure of that
class.
\end{Theorem}

\begin{Proof}
Consider the class $K$ of odd $4$-multipedes.  By Lemmas~4.3 and 4.4, $K$
is a finitely axiomatizable class of rigid structures.  By Lemma~4.7, for
every $\LL(C)$ sentence $\phi$, there exists $l$ such that $\phi$ does not
define a linear order in any $l$-meager $4$-multipede.  It remains to show
that $K$ contains an $l$-meager $4$-multipede.  By Theorem~3.1, there
exists an odd $l$-meager $4$-hypergraph $H$.  Extend $H$ to a $4$-multipede
by attaching two feet to each vertex of $H$, choosing positive slaves in
any way consistent with the definition of $2$-multipedes, ordering the
segments in an arbitrary way and finally adding representations of subsets
of segments.  The result is an $l$-meager $4$-multipede.
\end{Proof}

Call two structures $k$-equivalent if there is no $\Lk$ sentence which
distinguishes between them.  We answer negatively a question of Scott
Weinstein [W].

\begin{Theorem}
There exist $k$ and a structure $M$ such that every structure
$k$-equivalent to $M$ is rigid but not every structure $k$-equivalent to
$M$ is isomorphic to $M$.
\end{Theorem}

Theorem remains true even if $\Lk$ is replaced with $\Lk(C)$ in the
definition of $k$-equivalence. 

\begin{Proof}
By Lemma~4.3, there exists $k$ such that a first-order sentence with $k$
variables axiomatizes the class of odd $4$-multipedes.  By Theorem~3.1,
there exists a $2k$-meager odd hypergraph, and therefore there exists a
$2k$-meager odd $4$-multipede $M$.  By the choice of $k$, every structure
isomorphic to $M$ is rigid.  By Lemma~4.5, there a structure $k$-equivalent
to $M$ (even if counting quantifiers are allowed) but not isomorphic to
$M$.
\end{Proof}

\section*{References}

\begin{description}

\item[D]\ \ 
Anuj Dawar, ``Feasible Computation through Model Theory'', PhD Thesis,
Institute for Research in Cognitive Science, University of Pennsylvania,
Philadelphia, 1993.

\item[G]\ \ 
Yuri Gurevich, ``Logic and the challenge of computer science", In ``Current
Trends in Theoretical Computer Science" (Ed.\ E. B\"{o}rger), Computer
Science Press, 1988, 1--57.

\item[HKL]\ \ 
Lauri Hella, Phokion G. Kolaitis and Kerkko Luosto, ``How to Define a
Linear Order on Finite Models'', Symposium on Logic in Computer
Science, IEEE Computer Society Press, 1994, 40--49.

\item[IL]\ \ 
Neil Immerman and E.S. Lander, ``Describing Graphs: A First-Order Approach
to Graph Canonization'', in ``Complexity Theory Retrospective'', Ed.\ Alan
Selman, Springer Verlag, 1990, 59--81.

\item[S]\ \ 
Alexei Stolboushkin, ``Axiomatizable Classes of Finite Models and
Definability of Linear Order'', in Proc.\ 7th IEEE Annu.\ Symp.\ on Logic
in Computer Sci., Santa Cruz, CA, 1992, 64--70.

\item[W]
Scott Weinstein, Private Correspondence, Oct.\ 1993.

\end{description}

\end{document}